\documentclass[12pt]{article}

\usepackage{amsmath, amsfonts, amssymb}

\newtheorem{theorem}{Theorem}[section]
\newtheorem{lemma}[theorem]{Lemma}
\newtheorem{proposition}[theorem]{Proposition}
\newtheorem{corollary}[theorem]{Corollary}
\newtheorem{remark}[theorem]{Remark}

\newcommand{\CC}{{\mathbb C}}

\newcommand{\NN}{{\mathbb N}}
\newcommand{\ZZ}{{\mathbb Z}}

\newfam\msbfam
\font\tenmsb=msbm10 scaled \magstep1
\textfont\msbfam=\tenmsb
\font\sevenmsb=msbm7 scaled \magstep1
\scriptfont\msbfam=\sevenmsb
\font\fivemsb=msbm5 scaled \magstep1
\scriptscriptfont\msbfam=\fivemsb

\newcommand{\set}[1]{\{#1\}}

\newcommand{\calH}{{\mathcal H}}

\begin{document}

\title{
The Carath\'{e}odory-Cartan-Kaup-Wu theorem on an
infinite-dimensional Hilbert space}

\author{Joseph A. Cima, Ian Graham,\footnote{Author
supported in part by the Natural Sciences and Engineering Research
Council of Canada Grant \# A9221.} \\ Kang-Tae
Kim,\footnote{Author supported in part by the grant
KRF-2002-070-C00005 of The Korea Research Foundation.} and Steven
G. Krantz\footnote{Author supported in part by NSF Grant
DMS-9988854.}}

\maketitle

\begin{quote} \small
{\small \bf Abstract:} This paper treats a holomorphic
self-mapping $f: \Omega \rightarrow \Omega$ of a bounded domain
$\Omega$ in a separable Hilbert space ${\cal H}$ with a fixed
point $p$. In case the domain is convex, we prove an
infinite-dimensional version of the Cartan-Carath\'eodory-Kaup-Wu
Theorem. This is basically a rigidity result in the vein of the
uniqueness part of the classical Schwarz lemma. The main
technique, inspired by an old idea of H. Cartan, is iteration of
the mapping $f$ and its derivative. A normality result for
holomorphic mappings in the compact-weak-open topology, due to Kim
and Krantz, is used.
\end{quote}

\begin{quote}
\end{quote}

\setcounter{section}{-1}

\section{Introduction}

Perhaps the most important part of the classical Schwarz lemma is
the uniqueness statement: If $f: D \rightarrow D$ is a holomorphic
function from the unit disc $D$ to itself, $f(0) = 0$, and
$|f'(0)| = 1$, then $f$ is a rotation. This rigidity statement has
had considerable effect in the subject of complex differential
geometry. It is the wellspring of many holomorphically invariant
metrics, and has had a notable influence in the subject of mapping
theory.

The analogous result in higher (finite) dimensions has been
explored by Carath\'{e}odory, Cartan, Kaup, and Wu. See [KRA] and
[WU] for a careful discussion of the matter. The theorem is that,
if a holomorphic mapping $f:\Omega \to \Omega$ of a bounded domain
$\Omega$ in $\CC^n$ satisfies $f(p)=p$ for some $p \in \Omega$ and
$|\det (Df(p))|=1$, then $f$ is a biholomorphic mapping. We remark
here that the condition on the Jacobian determinant is, in this
case, equivalent to saying that every eigenvalue of the Jacobian
matrix has modulus 1.

In recent years, the study of holomorphic functions and mappings
on a Hilbert or Banach space has received much attention. Work of
Lempert (see for instance [LEM]) has served as a catalyst to this
activity. The more recent work of Kim/Krantz (see [KIK1], [KIK2])
is the more direct genesis of the present paper. In this rather
general setting, many of the familiar tools of finite-dimensional
analysis are no longer available. The geometry is much more
difficult. Yet there is interest in developing these ideas because
of potential applications to mathematical physics and partial
differential equations.

This paper develops a version of the
Carath\'{e}odory-Cartan-Kaup-Wu theorem on a separable Hilbert
space ${\cal H}$. We begin by formulating and proving a result on
a convex domain $\Omega \subseteq {\cal H}$.  This is Theorem 1.1.
Afterward, in Section 3, we present a more refined result
(Proposition 3.1) on the ball in a separable Hilbert space.

The authors thank the Banff International Research Station (BIRS)
for hosting us during the work on this problem, and for providing
a stimulating working environment. We thank John McCarthy, Peter
Rosenthal and Warren Wogen for advice and information about
analysis on Hilbert spaces.

\section{Statement of the main results}

Let $\calH$ be a separable, complex Hilbert space, and let $B$
denote its open unit ball. Let $S^1$ represent the unit circle in
the complex plane $\CC$. For a bounded linear operator $T$ on
$\calH$, denote by $\sigma(T)$ its spectrum, i.e.
$$
\sigma(T) = \set{ \lambda \in \CC \mid T - \lambda I
\hbox{ is not
invertible} } \, .
$$

A bounded linear operator $T: \calH \to \calH$ is said to be {\it
triangularizable} if there exists a basis $\set{e_1,e_2,\ldots}$
such that $T(\CC e_1 + \ldots + \CC e_N) \subset \CC e_1 + \ldots
+ \CC e_N$ for every positive integer $N$. In this case, we shall
sometimes say that $T$ is upper-triangular with respect to the
basis $\set{e_1, e_2, \ldots}$. Note that the basis that renders
an operator upper-triangular can always be taken to be
orthonormal, just because the Gram-Schmidt process preserves the
flags $E_N = \CC e_1 + \cdots + \CC e_N$.

Now we formulate the principal result of this paper:

\begin{theorem} \label{Main}  \sl
Let $\Omega \subseteq {\cal H}$ be a bounded convex domain.
Fix a
point $p \in \Omega$.  Let $f:\Omega \rightarrow \Omega$ be
a
holomorphic mapping such that
\begin{enumerate}
\item[{\bf (a)}]  $f(p) = p$ \, ; \item[{\bf (b)}]  the
differential $df_p$ is triangularizable; \item[{\bf (c)}]
$\sigma(df_p) \subseteq S^1$ \, .
\end{enumerate}
Then $f$ is a biholomorphic mapping.
\end{theorem}

To put the theorem in perspective, notice that in dimension $n =
1$ the hypothesis specializes to just $|f'(p)|=1$. In several (but
finitely many) variables, the hypothesis is equivalent to {\bf
(i)} non-degeneracy of the Jacobian matrix at $p$ and {\bf (ii)}
all eigenvalues of the Jacobian having modulus 1. In the
finite-dimensional case, one can always triangularize. In infinite
dimensions there are geometric conditions (involving the Fredholm
index) for triangularization.

Even if $T$ is a unitary operator on a Hilbert space ${\cal H}$
(so that certainly $\sigma(T) \subseteq S^1$), it does not
necessarily follow that $T$ is triangularizable. For example, take
$T$ to be the forward shift on $\ell^2(\ZZ)$. See [HER] for more
on the triangular operators.

Concerning the assumption {\bf (b)} in Theorem \ref{Main}, it
seems appropriate to present an example of a bounded convex domain
with an automorphism whose derivative at a fixed point is
upper-triangular.  Consider the map $A:\CC^2 \to \CC^2$ defined by
$A(z,w) = (z+bw, -w)$ for an arbitrary choice for $b \in
\mathbb{\CC} \setminus \set{0}$. Then $A$ is an involution (i.e.
$A \circ A = I$) with eigenvalues $\pm 1$; also $A$ is an
upper-triangular linear map. Take any bounded domain $D$ in
$\CC^2$ containing the origin. Let
$$
V = D \cup A(D).
$$
Then let $\Omega$ be the convex hull of $V$ in $\CC^2$. It is
obviously a bounded convex domain containing the origin. Then $A$
is an automorphism of $\Omega$ satisfying $A(0)=0$. Since $dA_0 =
A$, the differential $dA_0$ is upper-triangular.  Furthermore, if
one would sacrifice the convexity of the domain, it is possible
also to obtain a non-linear example. Let $G(z,w) = (z, w+z^2)$ for
instance, let $f = G \circ A \circ G^{-1}$, and let $W =
G(\Omega)$. Then $W$ is a bounded domain in $\CC^2$, and $f$ is an
automorphism of $W$ with $f(0)=0$. We see that $df_0$ is clearly
upper-triangularizable. It is simple to modify this example to
give an example in an infinite dimensional Hilbert space; one
simply needs to consider a map that is equal to $A$ on a
two-dimensional subspace and the identity on the orthogonal
complement of this subspace.

The referee has raised the question of how large is the class of
holomorphic mappings $f:\Omega\rightarrow \Omega$ satisfying the
conditions {\bf (a)}, {\bf (b)} and {\bf (c)} of Theorem 1.1,
within the class of such mappings satisfying only {\bf (a)} and
{\bf (c)}. We believe that rigidity phenomena for automorphisms of
domains in infinite dimensions would imply that the answer depends
very much on the domain. A related question is how large is the
class of triangularizable operators on a separable Hilbert space
$H$ within the class of all bounded linear mappings of $H$. In
this connection we mention the Weyl-von Neumann-Berg theorem [DAV,
p.\ 59], which asserts that every normal operator is a small
compact perturbation of a diagonalizable (in particular
triangularizable) operator. We thank the referee and the editors
for their remarks.

\section{Proof of Theorem~\ref{Main}}

The proof has several steps, some of which are of
independent interest.

\subsection{Basic Facts on the Differential}

Let $\Omega_1$ and $\Omega_2$ be domains in the separable Hilbert
space $\calH$ and let $f:\Omega_1\to\Omega_2$ be a holomorphic
mapping. The differential of $f$ at a point $p\in\Omega_1$ is
the bounded linear operator on $\calH$ defined by

$$
df_p (v) = \lim_{\CC\ni h \to 0} \frac{f(p+hv) -
f(p)}{h}.
$$

Write $B(p;r)=\{z\in \calH \mid \|z-p\| < r\}$.
Suppose that $\Omega_2$ is bounded, say $\Omega_2\subseteq B(0;M)$.
In this situation we have

\begin{lemma} \sl
If $p\in\Omega_1$ and $\hbox{dist}(p,\partial \Omega_1) = d >0$, then
$$\|df_p\|\le \frac{M}{d}\, .$$
\label{diffone}
\end{lemma}

\noindent\it Proof. \rm %
Choose $\rho > 0$ so that $\overline B(p;\rho) \subset \Omega_1$.
Let $v$ be a unit vector in $\calH$. The integral representation
\begin{equation} \label{Cauchyformula}
df_p(v)=\frac{1}{2\pi
i}\int_{|\zeta|=\rho}\frac{f(p+\zeta
v)}{\zeta^2} d\zeta
\end{equation}
leads immediately to the Cauchy estimates
$$
\|df_p(v)\|\le \frac{M}{\rho},
$$
and we may let $\rho$ tend to $d$.
\hfill $\Box\;$
\smallskip \\

\begin{lemma} \sl
If $\|df_z\|\le A$ for every $z\in B(p;r)$, then for any
pair $z,w\in B(p;r)$ one has
$$
\|f(z)-f(w)\| < 2rA \, .
$$
\label{difftwo}\end{lemma}

\noindent\it Proof. \rm %
Write
$$
f(z)-f(w)=\int_0^1 \frac{d}{dt}[f((1-t)w+tz)] dt
$$
and use the chain rule.
\hfill $\Box\:$
\smallskip \\

Now, in the situation of Theorem \ref{Main}, $f$ is a holomorphic
self-map of a bounded convex domain in $\calH$ and $p$ is a fixed
point. Let us write $T=df_p$.

We are assuming that $T$ is triangularizable. Hence it is
possible to choose a basis $e_1, e_2, \ldots$ such that
$$
T(\CC e_1 + \cdots + \CC e_N) \subset \CC e_1 + \cdots + \CC
e_N
$$
for every positive integer $N$. For convenience set
$$
E_N = \CC e_1 + \cdots + \CC e_N
$$
for every positive integer $N$. We shall call such $E_N$ a
\it
flag. \rm The union of these flags yields a vector space
that is
dense in $\calH$.
\smallskip \\

{From} here on we assume that the operator $T$ is
upper-triangular
with respect to a {\sl fixed} orthonormal basis $e_1,e_2,
\ldots$.
\smallskip \\

Notice that $\sigma(T)$ is contained in the unit circle. Since $T$
is upper-triangular, the diagonal entries are contained in its
spectrum $\sigma(T)$, and hence are of modulus 1.
\smallskip \\

We also note the following consequence of Lemma \ref{diffone}:

\begin{corollary} \sl
There exists a constant $C$ such that
$$
\|T^m\| \le C,
$$
for every positive integer $m$. \label{corol}
\end{corollary}

\subsection{Iteration of $T$}

In this section $T$ is an upper triangular matrix whose spectrum
$\sigma(T)$ is contained in $S^1$ and whose positive powers are
uniformly bounded in norm. Denote by $\lambda_j$ the $(j,j)$-th
diagonal entry of $T$. Then $|\lambda_j|=1$, $j\in\NN$.

\begin{lemma} \sl
There exists a sequence of natural numbers $\{m(k)\}_k$
such that for each fixed $j\in\NN$, $\lambda_j^{m(k)}\to 1$
as $k\to\infty$.
\label{T1}
\end{lemma}

\noindent \it Proof. \rm %
The sequence of powers $\{\lambda_1^k\}$ is bounded, so there
is a convergent subsequence, say $\lambda_1^{m(1,k)}\to \alpha_1$ as
$k\to\infty$, where $|\alpha_1|=1$. Similarly the sequence of powers
$\{\lambda_2^{m(1,k)}\}$ is bounded, so there is a subsequence
$\{m(2,k)\}$ of $\{m(1,k)\}$ such that
$\lambda_2^{m(2,k)}\to\alpha_2$ as $k\to\infty$, where $|\alpha_2|=1$.
Continuing in this way and using a diagonal sequence argument, we
obtain a subsequence $\{m(k,k)\}$ of the natural numbers and
complex numbers $\alpha_j$, $j\in\NN$, of modulus one such that
$\lambda_j^{m(k,k)}\to\alpha_j$ as $k\to\infty$ for each $j\in\NN$.
Now the sequence

$$m(k)=m(k+1,k+1)-m(k,k), \quad k\in\NN$$
is easily seen to have the property that
$$
\lambda_j^{m(k)}\to\frac{\alpha_j}{\alpha_j} = 1
$$
as $k\to\infty$, for each $j\in\NN$. \hfill $\Box\;$
\smallskip \\

If $A$ is an infinite square matrix we denote the $(i,j)$-th entry
by $A_{i,j}$, $1\le i,j\le\infty$. Also, we order the positions above
the diagonal in such a matrix first by column and then by row, i.e.
the ordering is $(1,2)$, $(1,3)$, $(2,3)$, $(1,4)\ldots$. If $N$
is a positive integer, we denote the $N\times N$ truncation of such
a matrix by $A_N$. The $N\times N$ identity matrix will be denoted
by $I_N$.

\begin{lemma} \sl
There is a sequence of natural numbers $\{\mu(k)\}_k$
such that, for each $N\in\NN$,
$$
(T_N)^{\mu(k)}\to I_N, \quad k\to\infty \, .
$$
\label{T2}
\end{lemma}

\noindent \it Proof. \rm %
Write $T = S+V$, where $S$ is a diagonal matrix with diagonal
entries of modulus 1 and $V$ has zeros on and below the main
diagonal. Denote the diagonal entries of $S$ by $\lambda_j,
j=1,2,\ldots.$

By Lemma \ref{T1}, there exists a sequence of natural
numbers $\{m(k)\}$ such that $\lambda_j^{m(k)}\to 1$ as
$k\to\infty$, for each fixed $j$.

The entries $\{(T^{m(k)})_{1,2}\}_k$ are bounded, so there is
a subsequence $\{\mu(k;1,2)\}$ of $\{m(k)\}$ such that
$\{(T^{\mu(k;1,2)})_{1,2}\}_k$ converges. By similar reasoning, a further
subsequence $\{\mu(k;1,3)\}$ has the property that
$\{(T^{\mu(k;1,3)})_{1,3}\}_k$ converges, and a further subsequence
of that one, denoted by $\{\mu(k;2,3\}$, has the property that
$\{(T^{\mu(k;2,3)})_{2,3}\}_k$ converges. Continuing in this way and
extracting a diagonal subsequence yields a subsequence $\{\mu(k)\}$
of $\{m(k)\}$ such that $\{(T^{\mu(k)}))_{i,j}\}_k$ converges for all
$(i,j),\quad 1\le i < j < \infty$.

Thus there is an infinite square matrix $W$, whose entries on and below
the main diagonal are zero, such that for each $N\in\mathbb{N}$
we have $(T_N)^{\mu(k)}\to I_N + W_N$.
(For fixed $N$ the convergence
may be taken to be in norm, but the norm convergence is not uniform in $N$.)

Now choose the smallest value of $N \ge 2$ with the
property that $W_{N}\ne 0$. Then at least one of the entries
in the last column of $W_N$ is nonzero, and all entries in the other
columns of $W_N$ are zero. If $\ell$ is a positive integer, then
$$((T_N)^{\mu(k)})^\ell = (T_N)^{\mu(k)\ell}
\to (I_N + W_N)^\ell, \quad k\to\infty,$$
and the entries in the last column of the matrix on the right
are given by $\ell$ times the corresponding entries in $(I_N + W_N)$
(except for a 1 in the $(N,N)$-position).

But this is a contradiction to the power boundedness of $T$ for
sufficiently large $\ell$ (see Corollary \ref{corol}).
We conclude that no such $N$ exists, i.e. $W=0$.
Therefore the sequence $\{\mu(k)\}$ has the property that
$(T_N)^{\mu(k)}\to I_N$ as $k\to\infty$ for each fixed $N$.
\hfill $\Box\;$

\subsection{Iteration of $f$}

Next, consider the iteration given by $f^1 = f$, $f^m =
f\circ
f^{m-1}$ for each integer $m > 1$.

We need the following two fundamental results. There is in fact a Banach
space version of the theorem of Kim and Krantz [KIK2]; we indicate
a proof here for the Hilbert space case.

\begin{theorem}[Kim/Krantz \hbox{ \bf [KIK2]}] \sl
Let $\Omega_1, \Omega_2$ be domains in a separable Hilbert space
$\calH$, and let $\Omega_2$ be bounded. Then every sequence
$\set{h_n : \Omega_1\to \Omega_2 \mid n=1,2,\ldots}$ of
holomorphic mappings admits a subsequence $\{h_{n(k)}\}_k$ that
converges to a holomorphic mapping $\widehat h$ from $\Omega_1$
into the closed convex hull of $\Omega_2$, in the
compact-weak-open topology (i.e., the compact-open topology in
which the strong topology is used on the domain space and the weak
topology is used on the target space). \label{KK}
\end{theorem}

\noindent \it Proof. \rm %
Let $\langle \, \cdot \, , \, \cdot \, \rangle$ denote the Hilbert
space inner product (linear in the first variable,
conjugate-linear in the second). One needs to show that for every
sequence $\set{h_n: \Omega_1 \to \Omega_2 \mid n=1,2, \ldots}$ of
holomorphic mappings, there is a subsequence $\{h_{n(k)}\}_k$ and
a holomorphic mapping $\widehat h : \Omega_1 \to \calH$ such that,
for all $g$ in the unit ball of $\calH$, $\langle
h_{n(k)},g\rangle \to \langle \widehat h,g \rangle$, uniformly on
compact subsets, as $k\to\infty$. If $g$ is such that $\Re \langle
w,g\rangle < 1$ for all $w\in\Omega_2$, then it is clear that $\Re
\langle \widehat h,g\rangle \le 1$, i.e. the image of $\widehat h$
must be contained in the closed convex hull of $\Omega_2$.

Let $\{z_n\}_{n\in\NN}$ be a dense sequence in $\Omega_1$.
We are going to do the usual diagonal sequence construction.
The sequence $\{h_n\}$ has a subsequence $\{h_{n(1,k)}\}_k$
such that $h_{n(1,k)}(z_1)$ converges weakly to an element
$\widehat h(z_1)$ as $k\to\infty$. This just uses the boundedness
of $\Omega_2$.

Now choose a subsequence of $\{n(1,k)\}_k$, denoted by
$\{n(2,k)\}_k$, so that $h_{n(2,k)}(z_2)$ converges weakly to
an element $\widehat h(z_2)$ as $k\to\infty$. Continue and then
choose the diagonal sequence $\{h_{n(k)}\}_k=\{h_{n(k,k)}\}_k$
generated by this process.

Since $\Omega_2$ is bounded, as we did earlier in the proof of
Lemma 2.1, we let $M$ be a positive constant such that $\Omega_2
\subset B(0; M)$.

Assume that $K$ is a compact subset of $\Omega_1$ and let
$\epsilon > 0$. Denote by $D= \hbox{dist }(K,\partial\Omega_1)$
and let $\delta =\min\{\frac{D}{3}, \frac{D\epsilon}{18M}\}$.
Cover $K$ with a finite number of balls $B(x_1;\delta),\ldots,
B(x_m; \delta)$ such that $B(x_\ell; \delta) \cap K \neq
\emptyset$ for any $\ell = 1,2,\ldots,m$, where each $x_\ell$
belongs to the dense sequence $\{z_n\}_{n\in\NN}$. Note that for
$z\in B(x_\ell;\delta)$ one has $\hbox{dist}(z,\partial
\Omega_1)
> D/3$. By Lemma \ref{diffone} we have $\|dh_n|_z \| <
\frac{3M}{D}$ for such $z$. This holds for every $n$. By Lemma
\ref{difftwo} and the choice of $\delta$, we have
\begin{equation} \label{equicontinuity}
\|h_n(z)-h_n(w)\| < \frac{6M\delta}{D} \le \frac{\epsilon}{3}
\end{equation}
for all $z,w\in B(x_\ell,\delta)$.
For any $g$ in the unit ball of $H$, choose $J=J(g)$ so that,
for $j,k > J$, we have
$$
|\langle h_{n(j)}(x_\ell)-h_{n(k)}(x_\ell),g\rangle| < \frac{\epsilon}{3}
$$
for $1\le \ell\le m$.

Now for any choice of $z\in K$ we have the existence of some $\ell$
with $z\in B(x_\ell,\delta)$.
Then for $j,k > J$ we have
\begin{eqnarray*}
|\langle h_{n(j)}(z)-h_{n(k)}(z),g\rangle|\, & \le & \,
|\langle h_{n(j)}(z)-h_{n(j)}(x_\ell),g\rangle| \\
&& \quad \, + \, |\langle h_{n(j)}(x_\ell)-h_{n(k)}(x_\ell),g\rangle|  \\
&& \quad \, + \, |\langle h_{n(k)}(x_\ell)-h_{n(k)}(z), g \rangle| \, .
\end{eqnarray*}
Each of the terms on the right hand side is less than
$\epsilon/3$.
For the first and third terms this follows from (\ref{equicontinuity}),
and for the second term it follows from the choice of $J$.

This shows that the sequence $h_{n(k)}(z)$ converges weakly to
the ``assignment'' $\widehat h(z)$ uniformly on compacta. It remains
only to show that the assignment is an analytic mapping on $\Omega_1$.
But given a fixed $z_0\in\Omega_1$ and a unit vector $v\in\calH$,
we can find $a>0$ such that the closed disc
$$
S=\{z_0+\zeta v \mid |\zeta|\le a\}
$$
is contained in $\Omega_1$, and hence the analytic functions
$\langle h_{n(k)}(z_0+\zeta v),g\rangle$ converge uniformly to
$\langle \widehat h(z_0+\zeta v),g\rangle$. The mapping
$\zeta\mapsto \widehat h(z_0+\zeta v)$ from the disc of radius $a$
in $\CC$ into $\calH$ is therefore holomorphic. This says
precisely that $\widehat h$ is Gateaux holomorphic [HIP], [MUJ].
Since a bounded Gateaux holomorphic mapping is holomorphic, we are
done.                                \hfill $\Box\;$

\begin{remark} \sl
By the Cauchy estimates, it follows from Theorem \ref{KK} that
the
derivative $dh_{n(k)}|_z (v)$ converges to $d\widehat h_z
(v)$
weakly, uniformly on compact subsets of $\Omega_1 \times
\calH$
(i.e., $z \in \Omega_1$ and $v \in {\cal H}$).
\end{remark}

To see this,
let $L$ be a compact subset of $\Omega_1\times\calH$.
Then for all $(z,v)\in L$,
we have $\hbox{dist}(z,\partial\Omega_1) \ge a > 0$ and $\|v\| \le b$,
for some positive constants $a$ and $b$.
It is elementary to see that there exists $r > 0$ such that
$$
K=\{z+\zeta v \mid (z,v)\in L, |\zeta|\le r \}
$$
is a compact subset of $\Omega_1$.
Now the relation (\ref{Cauchyformula}) gives
$$
dh_{n(k)}|_z(v) - dh_z(v) = \frac{1}{2\pi i}\int_{|\zeta|=r}
\frac{h_{n(k)}(z+\zeta v)-h(z+\zeta v)}{\zeta^2} d\zeta \quad ,
$$
for all $(z,v)\in L$. Hence for any linear functional $\tau$
on $\calH$,
$$
\tau\circ dh_{n(k)}|_z(v) - \tau\circ dh_z(v)
$$
$$
=\frac{1}{2\pi i}\int_{|\zeta|=r} \frac{\tau\circ h_{n(k)}(z+\zeta
v) -\tau\circ h(z+\zeta v)}{\zeta^2} d\zeta \, .
$$
The assertion follows easily from this in view of the compactness
of $K$.
\smallskip \\

\begin{theorem}[H. Cartan] \sl
Let $\Omega$ be a bounded domain in $\calH$ and let $p \in
\Omega$. Let $f:\Omega \to \Omega$ be a holomorphic mapping with
$f(p)=p$ and $df_p = I$. Then $f$ coincides with the identity
mapping of $\Omega$. \label{Cartan}
\end{theorem}

See [FRV], [KRA] for the proof of Cartan's theorem.
\smallskip \\

Now apply Theorem \ref{KK} to the sequence $\{f^{\mu(k)}\}$, where
$f$ is the mapping in Theorem 1.1 and $\{\mu(k)\}$ is the sequence
constructed in Lemma \ref{T2}. We obtain a subsequence
$\{f^{\nu(k)}\}$ that converges to some $\widehat{f}$ in the
compact-weak-open-topology. The sequence of derivatives at $p$ is
$\{T^{\nu(k)}\}$. The discussion of the preceding section shows
that
$$
\lim_{k\to\infty} (T_N)^{\nu(k)} = I_N
$$
for every positive integer $N$.

Thus
$$
d{\widehat f}_p|_{E_N} = I_{E_N} \hbox{ for every } N =
1,2,\ldots
\, .
$$
Note that $\displaystyle{\bigcup_N E_N}$ is dense in ${\cal
H}$.
Since $d\widehat{f}_p$ is bounded and $d\widehat{f}_p \bigr
|_{E_N} = I_N$ for all $N$, it follows that $d\widehat{f}_p
= I$.

Finally, the fact that $\widehat{f}(p) = p$ together with the
convexity of $\Omega$ implies that $\widehat f(\Omega)\subseteq
\Omega$. Now Cartan's theorem (Theorem \ref{Cartan}) implies that
$\widehat{f} \equiv \hbox{id}$.

\subsection{Proof of Theorem \ref{Main}}

It is time to complete the proof of Theorem \ref{Main}. From
the
last part of the preceding subsection, we have
\begin{equation} \label{tada-1}
\lim_{k \to \infty} f (f^{\nu(k)-1} (z)) = z = \lim_{k \to \infty}
f^{\nu(k)-1} (f(z))
\end{equation}
in the compact-weak-open topology on $\Omega$.
\smallskip \\

By Theorem \ref{KK}, choose a subsequence $\set{f_\ell}$ of
$\set{f^{\nu(k)-1}}$ that converges to a holomorphic mapping
$\widehat h:\Omega \to \Omega$ in the compact-weak-open topology.
(Recall that $\Omega$ is convex and $\widehat h(p)=p$).
\smallskip

The second identity in (\ref{tada-1}) implies that
\begin{equation}
\widehat h \circ f = \hbox{\rm id}. \label{tada-2}
\end{equation}
We see therefore that the holomorphic mapping $\widehat h$ is a
left inverse of $f$, and that $f \circ \widehat{h}$ is a
holomorphic mapping from $\Omega$ into itself.
\smallskip \\

On the other hand, one cannot immediately deduce from the first
identity in (\ref{tada-1}) that $f \circ \widehat h = \hbox{\rm
id}$, because it is only known at this point that $f_\ell$
converges to $\widehat h$ weakly.  (With respect to the weak
topology on the source-domain and the strong topology on the
target-domain, holomorphic mappings need not be continuous.) So it
is necessary to show that $f \circ \widehat h = \hbox{\rm id}$.
\smallskip \\

Now (\ref{tada-2}) implies that
$$
d{\widehat h}_p \circ df_p = I. \label{tada-3}
$$
Since $df_p$ is invertible, we see that $d{\widehat h}_p$ is
also
the right inverse of $df_p$. This implies that
$$
d(f \circ \widehat h)_p = df_p \circ d\widehat{h}_p = I.
$$
Applying Cartan's Theorem (Theorem \ref{Cartan}) again to
$f\circ
\widehat h:\Omega \to \Omega$, we see that
$$
f \circ \widehat h = \hbox{\rm id}.
$$
Therefore $f$ is a biholomorphic mapping of $\Omega$ onto itself.
This completes the proof. \hfill $\Box\;$

\section{Closing Remarks}

For holomorphic functions (i.e. $\CC$-valued functions) on a
domain in a separable Banach space, there is a normality result
for the compact-open topology, similar to the finite-dimensional
case [HUY], [KIK2], [MUJ]. For holomorphic mappings there is no
such result unless one makes further restrictions [HUY]. However,
it is possible to obtain interesting theorems about holomorphic
mappings using normality with respect to the compact-weak-open
topology.

The assumption that the differential be triangularizable at the
fixed point is not necessary for the conclusion of Theorem
\ref{Main} to be valid. In the case of the unit ball $B$ of the
Hilbert space $\calH$, a unitary map conjugated by a M\"obius
transformation is a holomorphic automorphism, but in general the
differential at the fixed point is not triangularizable.  This
phenomenon reflects the difficulty in the case of infinite
dimensional holomorphy caused by the excessive size of the
isotropy group in such cases as the ball (in finite dimensions,
large isotropy group characterizes the ball---see [GRK]). However,
the ball case has a reasonable formulation as follows.

\begin{proposition} \sl
Let $f:B \to B$ be a holomorphic mapping, let $p \in B$, and let
$M_p$ be a M\"obius transformation of the ball sending $p$ to the
origin. If $f$ satisfies:
\begin{itemize}
\item[\rm (1)] $f(p)=p$,
\item[\rm (2)] $U \circ d[M_p]|_p \circ df_p \circ d[M_p^{-1}]|_0$
is triangularizable, for some unitary transform $U$, and
\item[\rm (3)] $\sigma(U \circ d[M_p]|_p \circ df_p \circ
d[M_p^{-1}]|_0)$ lies in the unit circle in $\CC$,
\end{itemize}
then $f$ is an automorphism of $B$.
\end{proposition}

The arguments given in this article surely yield a proof of
Proposition 3.1. However, one can simplify the argument, thanks to
the fact that {\it the domain in consideration is the unit ball}.
Indeed, from Schwarz's lemma, upper triangularity and the spectrum
condition, one obtains that the operator $V = U \circ d[M_p]|_p
\circ df_p \circ d[M_p^{-1}]|_0$ has norm 1. It is known [DIN]
that such a linear transformation $V$ is in fact unitary. Thus
Cartan's Theorem applied to $V^{-1} \circ U \circ M_p \circ f
\circ M_p^{-1}$ implies that this map is the identity map. Hence
$f$ is a holomorphic automorphism of $B$. The authors would like
to thank Warren Wogen for pointing out this line of reasoning.

The convexity assumption on $\Omega$ was necessary due to the use
of weak convergence in several places.  Whether one can remove
this additional assumption should be an interesting problem to
explore in future work.  It would also be of interest to know
whether there is a result of Schwarz-Pick type in our
infinite-dimensional context.
\bigskip \bigskip \\

\noindent {\Large \bf References} \vspace*{.15in}

\begin{enumerate}

\item[{\bf [BOM]}] S. Bochner and W. T. Martin, {\it
Functions of
Several Complex Variables}, Princeton University Press,
Princeton,
1936.

\item[{\bf [CON]}] J. B. Conway, {\it A Course in Operator
Theory}, American Mathematical Society, Providence, RI,
2000.

\item[{\bf [DAV]}] K. R. Davidson, {\it $C^*$-algebras by
example}, Fields Institute Monographs, American Mathematical
Society, Providence, RI, 1996.

\item[{\bf [DIN]}] S. Dineen, {\it The Schwarz Lemma}, The
Clarendon Press, Oxford University Press, Oxford, 1989.

\item[{\bf [DUS]}] N. Dunford and J. T. Schwartz, {\it
Linear
Operators}, Interscience, New York, 1988.

\item[{\bf [FRV]}] T. Franzoni and E. Vesentini, {\it
Holomorphic
Maps and Invariant Distances}, North-Holland, Amsterdam,
1980.

\item[{\bf [GRK]}]  R. E. Greene and S. G. Krantz,
Characterization of complex manifolds by the isotropy subgroups of
their automorphism groups,  {\it Indiana Univ. Math. J.} 34(1985),
865-879.

\item[{\bf [HER]}] D.A. Herrero, Triangular operators, {\it Bull.
London Math. Soc.}, 23 (1991), 513-554.

\item[{\bf [HIP]}] E. Hille and R. S. Phillips, {\it Functional
Analysis and Semigroups}, Amer. Math. Soc. Coll. Publ. 31,
Providence, R. I., 1957.

\item[{\bf [HUY]}] C.-G. Hu and T.-H. Yue, Normal families
of holomorphic mappings, {\it J. Math. Anal. Appl.} 171(1992),
436-447.

\item[{\bf [KIK1]}] K.-T. Kim and S. G. Krantz,
Characerization of
the Hilbert ball by its automorphism group, {\it Trans.\
Amer.\
Math.\ Soc.} 354(2002), 2797--2818.

\item[{\bf [KIK2]}] K.-.T. Kim and S. G. Krantz, Normal
families
of holomorphic functions and mappings on a Banach space,
{\it
Expo.\ Math.} 21(2003), 193--218.

\item[{\bf [KRA]}] S. G. Krantz, {\it Function Theory of
Several
Complex Variables}, American Mathematical Society-Chelsea,
Providence, RI, 2001.

\item[{\bf [LEM]}] L. Lempert, The Dolbeault complex in
infinite
dimensions, {\it J. Amer.\ Math.\ Soc.} 11(1998), 485--520.

\item[{\bf [MUJ]}] J. Mujica, {\it Complex Analysis in Banach
Spaces}, North-Holland, Amsterdam and New York, 1986.

\item[{\bf [NAR]}] R. Narasimhan, {\it Several Complex
Variables},
University of Chicago Press, Chicago, 1971.

\item[{\bf [WU]}]  H. H. Wu, Normal families of holomorphic
mappings, {\it Acta Math.} 119(1967), 193--233.

\end{enumerate}
\vspace*{.2in}

\noindent Department of Mathematics \\
University of North Carolina \\
Chapel Hill, North Carolina 27514 \ \ USA  \\
{\tt cima@math.unc.edu}
\medskip \\

\noindent Department of Mathematics \\
University of Toronto \\
Toronto, CANADA \\
M5S 3G3 \\
{\tt graham@math.toronto.edu}
\medskip \\

\noindent Department of Mathematics \\
Pohang University of Science and Technology \\
Pohang 790-784 KOREA  \\
{\tt kimkt@postech.ac.kr}
\medskip \\

\noindent Department of Mathematics \\
Campus Box 1146  \\
Washington University in St. Louis \\
St.\ Louis, Missouri 63130 \ \ USA   \\
{\tt sk@math.wustl.edu}

\end{document}